# Application of support vector machine for the fast and accurate reconstruction of nanostructures in optical scatterometry

Jinlong Zhu*, Hao Jiang, Chuanwei Zhang, Xiuguo Chen, and Shiyuan Liu*

*State Key Laboratory of Digital Manufacturing Equipment and Technology,
Huazhong University of Science and Technology, Wuhan 430074, China*

*Corresponding author: Leon.J.Zhu2@gmail.com; shyliu@hust.edu.cn

Nonlinear regression methods such as local optimization algorithms are widely used in the extraction of nanostructure profile parameters in optical scatterometry. The success of local optimization algorithms heavily relies on the estimated initial solution. If the initial solution is not appropriately selected, it will either take a long time to converge to the global solution or will result in a local one. Thus, it is of great importance to develop a method to guarantee the capture of global optimal solution. In this paper, we propose a method combining support vector machine (SVM) and Levenberg-Marquardt (LM) algorithm for the fast and accurate parameters extraction. The SVM technique is introduced to pick out a sub-range in the parameters' rough range, in which an arbitrary selected initial solution for the LM algorithm is then able to achieve the global solution with a higher possibility. Simulations and experiments conducted on a one-dimensional Si grating and a deep-etched multilayer grating have demonstrated the feasibility and efficiency of the proposed method.

## 1. INTRODUCTION

Scatterometry is one of the most promising metrology techniques for semiconductor manufacturing industry because it meets the requirements of high-volume manufacturing and in-line monitoring when compared with traditional methods such as scanning electron microscopy (SEM) and atomic force microscopy (AFM) [1, 2]. Generally, this technique involves two main procedures, namely the forward modeling and the inverse problem solution. For the forward modeling, several techniques such as rigorous coupled-wave analysis (RCWA) [3-5], the boundary element method (BEM) [6], and the finite-difference time-domain (FDTD) method [7] have been introduced to simulate the scattering signature from a nanostructure. Here the general term signature contains the scattered light information of the diffractive structure, which can be in the form of reflectance, ellisometric angles, Stokes elements, and Mueller matrix elements. The second procedure involves the extraction of profile parameters from the measured signature, which is a typical inverse problem with the objective of finding the profile parameters whose calculated signature can best match the measured one.

To solve the inverse problem in scatterometry, generally, two kinds of approaches are widely used, i.e., the library search [8-11] and the nonlinear regression [12-14]. The library search is also called the look-up table, which is based on the building and retrieval of a signature library. To know more about the library search technique, we recommend the readers to [8]. The nonlinear regression can be further divided into two types of optimization methods, namely, the global optimization methods (GOM) such as simulated annealing [15] and genetic algorithm (GA) [16], and the local optimization methods (LOM) such as Gauss-Newton and Levenberg-Marquardt (LM) algorithms [17]. The GOM can guarantee the global solution theoretically, but the time consumption is usually unacceptable. Therefore the LOM is more commonly used in scatterometry due to its fast convergence rate. However, the LOM depends heavily on the initialization, that is, an inappropriate selection of the initial estimation of profile parameters will either result in more calculation time or will lead to a local solution. Traditionally the issue of initialization is dealt with by broadening the range of wavelength and using more wavelength points in a signature [18]. However, the increase of wavelength points in a signature not only accompanies with the penalty of computational cost, but also brings other problems such as "how many wavelength points in a signature are sufficient to ensure the global solution". As an alternative, we may use some artificial intelligence methods to give a reliable initial estimation of profile parameters. Zhang et al. proposed a method combining the ANN and LM algorithm for robust and fast extraction of profile parameters of deep trench structures, in which the ANN is used to estimate an initial value [19]. Gereige et al. proposed a method based on ANN classification to give a first estimation of the solution location in the multi-parameter space, and this method has demonstrated that it can efficiently offer a reliable and qualitative characterization of the sample in term of layer thickness [20]. In addition to the ANN, the support vector machine (SVM) [21-23] is another kind of artificial intelligence technology, which is based on the statistical learning theory [24, 25], and displays a better generalization performance than ANN [21, 22]. SVM is a technique developed for the task of classification, and was originally designed to solve the binary classification problems. Its basic idea is to use a trained support vector network to non-linearly map the input vectors to a very high-dimension feature space, in which a linear decision surface can be constructed to divide the mapped input vectors into two groups. In consideration of the fact that most of the classification problems in practice can be attributed to a multi-classification one, researchers have developed several multi-classification SVM algorithms such as "one-against-all", "one-against-one", and directed acyclic SVM [26]. The developed multi-classification ability of SVM largely expands its application range. Jin et al. used the SVM to approximate the relation between scattering signature and nanostructure profile parameters [27]. Baek et al. used the SVM to determine whether a semiconductor sample should or should not proceed to subsequent manufacturing steps according to the mapping results of diffraction signals [28]. We recently applied the SVM-based classification to accurately recognize the geometrical profile of a one-dimensional (1D) trapezoidal grating and developed an SVM-based library search strategy [9]. In the present work, we apply the SVM technique to deal with the issue of initialization. The pre-trained SVM classifiers are used to pick out a sub-range from the parameters' rough range, and then an arbitrary selection of initial solution in the sub-range for the LM algorithm is able to ensure the global solution with a relative high possibility.

The reminder of this paper is organized as follows. Section 2 briefly describes the principle of SVM and the combined method. Section 3 presents the results, including the measurement setup and sample description, numerical and experimental results. Finally, we draw some conclusions in Section 4.

## 2. METHOD

### 2.1 Principle of SVM

SVM was firstly introduced by Vapnik [24, 25] to solve the binary classification problem. After that the SVM has undergone a great development. We begin with the problem of binary classification for the purpose of describing the principle of SVM. For a binary classification problem, the training pairs are represented as

$$(\mathbf{x}_1, y_1), ..., (\mathbf{x}_i, y_i), ..., (\mathbf{x}_N, y_N), \mathbf{x}_i \in R^n, y_i \in \{-1, 1\}, i = 1, 2, ..., N, \quad (1)$$

where $\mathbf{x}_i$ is an $n$-dimensional vector representing the $i$th training signature, $y_i$ is a scalar with two values of -1 and 1 representing two classes. The combination is called a training pair, and $N$ is the number of training pairs. The training pairs are used in the training of an SVM classifier.

For a measured signature $\mathbf{x}$, the value of a decision function $f(\mathbf{x})$ determines which class $\mathbf{x}$ belongs to. The decision function can be expressed as

$$f(\mathbf{x}) = \text{sign}\left[\mathbf{w} \cdot \boldsymbol{\psi}(\mathbf{x}) + b\right], \quad (2)$$

where $\boldsymbol{\psi}(\mathbf{x})$ is a mapping function of $\mathbf{x}$, $b$ is a bias, and $\mathbf{w}$ is a support vector that can be expressed as a linear combination of $\boldsymbol{\psi}(\mathbf{x}_i)$:

$$\mathbf{w} = \sum_{i=1}^{N} \lambda_i y_i \boldsymbol{\psi}(\mathbf{x}_i), \quad (3)$$

where $\lambda_i$ is the weight coefficient of the $i$th input training signature. By substituting Eq. (3) into Eq. (2), and by defining a new function

$$k(\mathbf{x}, \mathbf{x}_i) = \boldsymbol{\psi}(\mathbf{x}_i) \cdot \boldsymbol{\psi}(\mathbf{x}), \quad (4)$$

we can get the final expression of the decision function

$$f(\mathbf{x}) = \text{sign}[\sum_{i=1}^{N} \lambda_i y_i k(\mathbf{x}, \mathbf{x}_i) + b]. \quad (5)$$

The function $k(\mathbf{x}, \mathbf{x}_i)$ in Eqs. (4) and (5) is called the kernel function, which plays an important role in SVM. We list several kernel functions that will be studied in this paper. The first one, polynomial kernel which is expressed as

$$k_p(\mathbf{x}, \mathbf{x}_i) = (\mathbf{x} \bullet \mathbf{x}_i)^d, \quad (6)$$

where $d$ is the power exponent. Compared to the polynomial kernel, the radial basis function (RBF) kernel and Sigmoid kernel are two kernels that are more nonlinear and expressed as

$$k_r(\mathbf{x}, \mathbf{x}_i) = \exp(-\frac{\|\mathbf{x} - \mathbf{x}_i\|_2}{\sigma}), \quad (7)$$

$$k_s(\mathbf{x}, \mathbf{x}_i) = \tanh(\beta \mathbf{x} \cdot \mathbf{x}_i^T), \quad (8)$$

respectively. The $\sigma$ in Eq. (7) and $\beta$ in Eq. (8) respectively are the scaling factors, T represents the transpose. In this paper, for simplicity, the $d$, $\sigma$ and $\beta$ are called controlling factors.

In consideration of multi-classification problem discussed in this article, for simplicity, we directly use the advanced support vector machines tool introduced by Chang and Lin for multi-classification [29].

### 2.2 The SVM/LM Combined Method

Supposing the profile of a nanostructure corresponds to $M$ profile parameters, and each parameter has a priori rough range $R_i$ ($i = 1, 2, …, M$). Each priori rough range $R_i$ is sliced into $N_i$ ($i = 1, 2, …, M$) sub-ranges with each sub-range represented by a unique symbol. We then take rough range $R_i$ of $i$th parameter for example: $R_i$ is sliced into $N_i$ sub-ranges, and we randomly generate $K_j$ ($j = 1, 2, …, N_i$) points in each sub-range. Also, for the rest $M$-1 rough ranges of parameters, we randomly generate $K_2$ - $K_M$ points respectively. Hence for an arbitrary sub-range $K_j$ of $R_i$ and the rough ranges of the other $M$-1parameters, we can make $K_j \prod_{i=2}^{M} K_i$ combinations in total and each combination completely characterizes the profile of the nanostructure. Our in-house RCWA software package is used to calculate the corresponding $K_j \prod_{i=2}^{M} K_i$ optical signatures for the $K_j \prod_{i=2}^{M} K_i$ combinations, and these $K_j \prod_{i=2}^{M} K_i$ signatures with their unique symbol form a training set. Since the range of $i$th parameter $R_i$ has been sliced into $N_i$ sub-ranges, we can generate $N_i$ training sets for each parameter. Hence we can use the $N_i$ training sets to train an SVM classifier with $N_i$ classes, and each class corresponds to a training set. The trained SVM classifier can map the measured signature of a nanostructure into its corresponding sub-range of $i$th parameter. Therefore, we can generate all the $M$ SVM classifiers for the $M$ parameters. Once the measured signature of a nanostructure is inputted into the $M$ trained SVM classifiers respectively, we can obtain $M$ symbols with each symbol corresponds to a sub-range of each parameter. Each mapped symbol indicates that the true value of a profile parameter locates in the corresponding sub-range with a relative high possibility. Since the sub-range is smaller than the corresponding rough range, it contains less local solutions. As such, an initial solution selected in these sub-ranges for the LM algorithm is able to achieve the global solution with a relative high possibility. The flowchart of profile parameter extraction with the SVM and LM combined method is presented in Fig. 1.

## 3. RESULTS

### 3.1 Measurement Setup and Sample Description

A dual-rotating compensator Mueller matrix ellipsometer (DRC-MME) suitable from ultraviolet to infrared spectrum (ME-L ellipsometer, Wuhan Eoptics Co., Ltd.) is used for demonstration. Data analysis is performed using the in-house developed optical modeling software based on rigorous coupled-wave analysis (RCWA) [3-5]. As shown in Fig. 2, the system configuration of the DRC-MME in order of light propagation is PCr₁SCr₂A, where P and A stand for the fixed polarizer and analyzer, $Cr_1$ and $Cr_2$ refer to the first and second frequency-coupled rotating compensators, and S stands for the sample [30, 31]. With the light source used in this ellipsometer, the wavelengths available are in the 200 and 1000 nm range, covering the spectral range of 200 to 800 nm used in this article. The azimuth angle and incident angle in this article are chosen as 0° and 65° respectively.

Two samples are studied for demonstration. The first sample is a (100)-orientation trapezoidal Si grating, whose cross-section image obtained by an SEM (X-SEM) (Nova NanoSEM450, FEI Co.) is presented in Fig. 3(a). The grating parameters measured by SEM are the top critical dimension $TCD$ = 350 nm, the height of grating $Hgt$ = 472 nm, the bottom critical dimension $BCD$ = 383 nm, and $Pitch$ = 800 nm. In the following simulation and experiment parts, the period Pitch is fixed at 800 nm for simplicity. The etched Si grating is chosen for this study due to its long-term dimensional stability, high refractive index contrast, and relevance to the semiconductor industry [32]. Optical properties of Si are taken from [33].

The second sample is a deep-etched multilayer grating, whose cross-section image measured by a transmission electron microscopy (TEM) (TE20, TEM.FEI Co.) is presented in Fig. 3(b). This sample consists of Si,

SiO$_2$, and Si$_3$N$_4$ trapezoidal gratings from bottom to top. The overall parameters under measurement include the top critical dimensions (CD) $D_1$ and $D_2$ of the Si$_3$N$_4$ layer and SiO$_2$ layer respectively, the bottom CD $D_3$ of the Si layer, and the thicknesses $H_1$, $H_2$, $H_3$ of Si$_3$N$_4$ layer, SiO$_2$ layer, and Si layer, respectively. The TEM measured values of the above parameters are $D_1$ = 75.01 nm, $H_1$ = 135.60 nm, $D_2$ = 86.90 nm, $H_2$ = 9.92 nm, $D_3$ = 124.13 nm, $H_3$ = 134.29 nm, and pitch $P$ = 154 nm.

### 3.2 Numerical Results

The effectiveness of the proposed method combing SVM and LM algorithm mainly depends on the classification ability of the SVM classifiers. In practice, many factors, such as the number of training pairs in each training set, the number of wavelength or incident angle points in one training signature and different kernels with different controlling factors, may affect the classification ability of the trained SVM classifiers. Thus it is of great importance to evaluate the effects of all these factors. Here we first evaluate the effects of wavelength point number in a training signature and kernel selections on the classification accuracy of SVM classifiers. The classification accuracy is defined as the ratio between the number of the correct classified signatures and total number of the input testing signatures.

Now we suppose $TCD$, $Hgt$ and $BCD$ of the 1D Si grating are the parameters that need to be extracted, which means $M$ = 3, and their rough ranges are set as 250 - 550 nm, 300 - 600 nm and 250 - 550 nm respectively. We then set $N_1$, $N_2$ and $N_3$ as 4, i.e., each rough range is equally divided into four sub-ranges. We take the training of the first SVM classifier which corresponds to parameter $TCD$ for example. The rough range of $TCD$ is divided into 4 sub-ranges with equal lengths, which means the first SVM classifier has 4 classes. For the generation of the training set for the 1st class, we randomly generate 15 points in the first sub-range of $TCD$ ($K_1$ = 15), and randomly generate 15 points in each rough range of the other parameters respectively ($K_2$=$K_3$=15). Thus we can obtain $\prod_{i=1}^{3} K_i$ = 3375 combinations in total, and then use our in-house RCWA software package to generate the corresponding 3375 Mueller matrices to form the training set. After that we then use the same procedure above to generate 3375 training pairs for each of the rest three classes respectively.

Once the four training sets have been prepared we can train the first SVM classifier for parameter $TCD$. The same procedure is used to train the other two SVM classifiers that correspond to parameters $Hgt$ and $BCD$, respectively. Once the three SVM classifiers have been trained successfully, one hundred testing Mueller matrices that are randomly generated within the rough ranges of $TCD$, $Hgt$ and $BCD$ are input into the SVM classifiers to evaluate the classification accuracies under different numbers of wavelength points and different kernels with different controlling factors. The wavelength points are equally selected in the wavelength range 200 - 800 nm. The results are presented in Fig. 4.

In Fig. 4(a), we can observe that, for the case of $d$ = 1 in polynomial kernel, the classification accuracy is arising with the increase of wavelength point number. Then, with the increase of $d$ (up to 4), the classification accuracies correspond to all the different numbers of wavelength points are arising. Further, if $d$ reaches 5, the classification accuracy corresponds to 7 wavelength points keeps arising and reaches 98%, while for the other different numbers of wavelength points their corresponding classification accuracies slightly decrease. The results present in Fig. 4(a) have demonstrated that the number of wavelength points is not necessary to be as large as possible to achieve the good classification accuracy, and sometimes the large number of wavelength points may just be the opposite to what one wishes. This phenomenon is due to the SVM overfitting, which arises from the fact that only a small number of wavelength points have meaningful contributions to data variation [34, 35]. Hence, we should emphasize, when polynomial kernel is used and $d$ is set as 5, a training Mueller matrix with only 7 wavelength points is enough to achieve the relatively high classification accuracy (98%). This is a significant advantage since the reduction of wavelength point number will largely reduce the time consumption. If RBF kernel is selected, as shown in Fig. 4(b), only when $\sigma$ = 1 we can obtain the relatively high classification accuracies for all the different numbers of wavelength points. Again, we find that the large number of wavelength points does not ensure the high classification accuracy, and even the least number of wavelength points (7 wavelength points) can ensure the high classification accuracy (94%) at the condition of $\sigma$ = 1. While for the Sigmoid kernel, as shown in Fig. 4(c), in the range from 0.01 to 100 of the controlling factor $\beta$, all of the classification accuracies are below 60%, which indicates, at least in the range of $\beta$ shown in Fig. 4(c), the Sigmoid kernel is not suitable for the classification cases in this paper. Thus, we will not take the Sigmoid kernel into consideration in the following content.

The above simulation has presented the influence of different numbers of wavelength points and different kernel functions with different controlling factors on the classification accuracy. In the present paragraph we will present the relationship between the number of training pairs in a training set and the classification accuracy. As have discussed in the above paragraph, the classification accuracy can be very high even the number of wavelength points of a training signature is as small as 7, on the condition of the appropriate selected controlling factors. Hence, for the simulations here we set the number of wavelength points as 7. Also, we set $d$ = 5, $\sigma$ = 1 for polynomial and RBF kernels respectively. The simulation results are presented in Fig. 5. It is obvious that with the increase of training pair number the classification accuracy gets higher and higher. In summary, we can point out that for the polynomial kernel the 7 wavelength points is able to ensure 98% classification accuracy if there are 3375 training pairs in a training set and $d$ is set as 5. Moreover, for the RBF kernel the 7 wavelength points is still able to ensure the 94% classification accuracy if there are the same number of training pairs as above and $\sigma$ is set as 1.

The above simulation is conducted in the ideal condition. However, in practice the measured signature suffers from the random errors and offset errors [36], which may finally affect the classification accuracy. Hence, it is of great importance to evaluate the effect of different measurement errors on the classification accuracy. The generation approach of random errors and offset errors can be found in [36], in which the standard deviations of the simulated random errors and offset error values at different wavelengths are set as several percents of root-mean-square (rms) in the Mueller matrix over the full wavelength range of interest. In this paper, we use the term "magnitude" instead of the "several percents of rms" to describe the strength of random error and offset error. The effects of random errors and offset errors on the classification accuracy are shown in Fig. 6. It is observed that the RBF kernel presents the better resistance to both the random errors and the offset errors than the polynomial kernel.

The above simulations have demonstrated that the SVM classifiers obtained by using 3375 training pairs with 7 wavelength points in each training signature and the RBF kernel with $\sigma$ = 1 is capable to achieve the relative higher classification accuracy. Then, we will further demonstrate the effectiveness of the proposed method which combines SVM and LM algorithm.

In order to test the capability of the proposed SVM/LM combined method for different nanostructures, a data set of 100 "measured" Mueller matrix was generated. To generate the 100 "measured" Mueller matrix, 100 combinations of profile parameters $TCD$, $Hgt$ and $BCD$ were generated randomly and independently in the rough ranges of 250 - 550 nm, 300 - 600 nm and 250 - 550 nm respectively. Then the corresponding 100 Mueller matrices of the 100 nanostructures were generated using RCWA in the wavelength range 200 - 800 nm with an increment 10 nm. We added random errors and offset errors whose magnitudes are both 0.05 into the corresponding generated Mueller matrices to simulate the "measured" signatures. We then picked the signature wavelength points at wavelengths 200 nm, 300 nm, 400 nm, 500 nm, 600 nm, 700 nm and 800 nm to form the to be mapped signature for each Mueller matrix and input it into the SVM classifiers. As expected, the three trained SVM classifiers mapped the "measured" Mueller matrices into the corresponding sub-ranges of *TCD*, *Hgt* and *BCD* with the classification accuracies 97%, 99% and 100% respectively. After the mapping, the medians of the mapped sub-ranges of *TCD*, *Hgt* and *BCD* for each "measured" Mueller matrix then were taken as the initialization for the LM algorithm. To make a comparison, we also performed the directly LM iterations by using the medians of the three parameter rough ranges as the initialization for the 100 "measured" Mueller matrices. The time consumptions together with the errors of extracted *TCD*, *Hgt* and *BCD* by the SVM/LM combined method and the LM directly are presented in Fig. 7(a) and Fig. 7(b) respectively. As expected, in Fig. 7(a) the time consumptions for the 100 "measured" Mueller matrices by the SVM/LM combined method are all about 100 s for the 100 "measured" Mueller matrices obtained on the 2.3 GHz Intel i5-2410M personal computer. Moreover, the errors of extracted *TCD*, *Hgt* and *BCD* by the SVM/LM combined method are all within the range 0 - 2 nm, which implies a high "measurement" accuracy. The high resolution indicates that the SVM classifiers have accurately mapped the 100 "measured" Mueller matrices into the sub-ranges which contain the global solution, and choosing the medians of sub-ranges as the initialization is reasonable. In contrast, in Fig. 7(b) some of the time consumptions by LM algorithm directly are as large as 1250 s, and they usually correspond to the large errors of extracted *TCD*, *Hgt* and *BCD*. The large "measurement" errors indicate that the initialization for the corresponding "measured" Mueller matrices is not the appropriate one, thus the LM iteration results are non-convergence. The above simulations have demonstrated that the SVM mapping is able to provide an appropriate initialization for the LM algorithm, which helps avoid the non-convergence with a high possibility. The issue of initialization is crucial, since the convergence domain of global solutions are smaller for most cases provided the number and rough ranges of profile parameters are larger. To further demonstrate the superiority of the SVM/LM combined method, we pick out those convergent results by the LM algorithm among the 100 testing results and compare them with that of the SVM/LM combined method. As can be seen in Fig. 8, obviously, all the extracted errors of *TCD*, *Hgt* and *BCD* of SVM/LM are the same with those of the LM, which are expected by the deterministic property of the LM algorithm. However, most of the time consumptions of the SVM/LM method are less than those of the LM, indicating that the initializations obtained by the SVM mapping are closer to the global solutions than the medians of those rough ranges. The less time consumption can be achieved by pre-dividing the rough ranges of parameters into more sub-ranges.

We next consider the simulation for the second sample, which is a deep-etched multilayer grating with cross-section image shown in Fig. 3(b). The true values of the geometrical parameters are respectively set as $D_1 = 75$ nm, $H_1 = 135$ nm, $D_2 = 86$ nm, $H_2 = 10$ nm, $D_3 = 124$ nm and $H_3 = 134$ nm. Following the same process and the established experience we choose the RBF kernel with $\sigma = 1$ for the training of six SVM classifiers (note that six geometrical parameters are used to fully characterize the sample). The rough ranges of the geometrical parameters: $D_1 \in [55$ nm, $95$ nm$]$, $H_1 \in [115$ nm, $155$ nm$]$, $D_2 \in [66$ nm, $106$ nm$]$, $H_2 \in [1$ nm, $19$ nm$]$, $D_3 \in [104$ nm, $144$ nm$]$ and $H_3 \in [114$ nm, $154$ nm$]$, are equally divided into four segments, which means the corresponding SVM classifiers have 4 classes. $10^6$ training pairs are generated using the same in-house RCWA software package to form the training set for each class of an SVM classifier. Here 21 points are equally selected in the wavelength range of 200 - 800 nm for the tradeoff of both overdetermination and computation efficiency. Once the six SVM classifiers have been generated, one hundred "contaminated" Mueller matrices that are randomly generated within the rough ranges of all the geometrical parameters are input into the SVM classifiers to evaluate the classification accuracy, as presented in Fig. 9 (a). Obviously, with the increase of noise level and offset magnitude, the classification accuracy presents a decreasing trend, but even in the worst case the classification accuracy value still stays above 84%. We next compare the computation time for the SVM/LM and LM methods. For the SVM/LM combined method, a Mueller matrix is firstly inputted into the SVM classifiers to get the corresponding sub-ranges, after which the medians of these sub-ranges are set as the initialization for the rest LM iteration. While for the sole LM iteration the left margins of the rough ranges of geometrical parameters are set as the initialization. In the reconstruction process the wavelength range is still kept as 200 - 800 nm but the increment is reset as 10 nm to strengthen the iterative robustness. It can be observed in Fig. 9(b) that the SVM/LM combined method takes about 102.4 s to convergent to the preset tolerance, while all the time consumptions of LM method are beyond 102.8 s, and meanwhile portions of the results are even larger than 103.4, which indicates the divergent solutions. The above simulation establishes the feasibility of the proposed SVM/LM combined method in the fast and accurate reconstruction of complex nanostructure.

### 3.3 Experimental Results

The measured Mueller matrices corresponding to the 1D Si grating and the deep-etched multilayer grating are both obtained in the wavelength points from 200-800 nm with 10 nm increment under 0° azimuth angle and 65° incident angle by the DRC-MME. We then picked the wavelength points which correspond to the 7 spectral points (200 - 800 nm with 100 nm increment), and used the optimal SVM classifiers trained under the conditions of 3375 training pairs with 7 wavelength points and RBF kernel with $\sigma = 1$ and 4 classes, to map the measured Mueller matrix of the 1D grating. Besides, we also use the six SVM classifiers generated in the above section to map the measured Mueller matrix of the deep-etched multilayer grating. The mapping results show that the true geometrical values of 1D Si grating are within the sub-ranges which are respectively 325 - 400 nm, 450 - 525 nm and 325 - 400 nm, and the true values of the deep-etched multilayer grating are respectively within the subranges 75 - 85 nm, 125 - 135 nm, 86 - 96 nm, 14.5 - 19 nm, 124 - 134 nm and 124 - 134 nm. We then chose the medians of the mapped sub-ranges as the initializations and input them into the LM algorithm. The extracted profile parameters together with the time consumption of the 1D Si and deep-etched multilayer gratings are presented in Table 1 and Table 2, respectively.

In the first test case of 1D Si grating, four additional parameter sets are selected as the initializations, which are (400, 450, 400) nm, (520, 430, 450) nm, (500, 525, 500) nm and (475, 525, 500) nm respectively, for comparison. The initialization (400, 450, 400) nm typically consists of the medians of the rough ranges of *TCD*, *Hgt* and *BCD*. The four initializations were inputted into the LM respectively, and their corresponding results are also presented in Table 1. We can find that the extracted profile parameters using SVM/LM combined method is the same as the LM algorithm with the initialization (400, 450, 400) nm, and agree to the results reported by SEM. The extracted profile parameters achieved by the LM algorithm with the rest initializations such as (520, 430, 450) nm, (500, 525, 500) nm and (475, 525, 500) nm show obvious deviations from the SEM measurements, which indicates the LM iterations fall into incorrect local solutions. Moreover, the time consumption of the SVM/LM with initialization (362.5, 487.5, 362.5) nm is only 122.5 s, which is significantly less than the time consumptions of the rest four cases. These results show that the inappropriately selected initial values for LM method lead to the local solutions, while the proposed method is capable to avoid such a situation and faster the convergence by limiting the initial values in a small range determined by SVM. Further, we calculated the Mueller matrices of the extracted profile parameters in the first and third rows of Table 1 respectively, and presented them in Fig. 10, in which we find the Mueller matrix corresponds to the SVM/LM method fits the measured Mueller matrix very well. The slightly misfit of the off-diagonal elements may be due to the inevitable random error and actual nonzero azimuthal angle, since in our DRC-MME prototype there is no positioning device on the rotating platform yet and the zero azimuthal angle is roughly guaranteed by the manual mode. Actually, the misfit of the off-diagonal elements did further demonstrated the SVM classifiers are robust to the system error and random error since the training of SVM classifiers only takes advantage of the simulated signatures. However, for the one corresponds to the LM algorithm with initialization (520, 430, 450) nm, the calculated Mueller matrix biases the measured Mueller matrix largely.

In the second test case of deep-etched multilayer grating, we can find that the SVM/LM combined method leads to the convergent reasonable solution within 354.78 s, while an inappropriate initialization is ease to result in not only the divergent solution but also the consumption of huge time, as shown in Table 2. Slight changes of few elements in the initializations (95, 155, 106, 19, 144, 154) nm and (55, 115, 66, 1, 104, 114) nm result in divergent results, which significantly deviate from the TEM measured ones. The results we achieved not only demonstrates the feasibility of the proposed SVM/LM method, but also show us the initialization selection may become a more critical factor in the profile reconstruction for complex nanostructures.

The above experimental results have demonstrated that the SVM classifiers are able to correctly classify the measured signature of a nanostructure with unknown parameter values into the sub-ranges of the rough ranges, which, reduces the possibility of suffering from local solutions. Moreover, the mapped sub-ranges indeed shorten the solution space, in which an arbitrary chosen initialization for the LM algorithm is more possibly able to converge to the global solution with higher speed than the directly using LM in the rough ranges.

## 4. CONCLUSION

In this paper, in order to deal with the issue of initialization in the solution to the inverse problem in optical scatterometry, we have introduced the SVM classifiers to map the measured Mueller matrices of a 1D etched Si grating and a deep-etched multilayer grating to the pre-divided sub-ranges of the rough ranges. Then the initializations selected in the mapped sub-ranges for the LM algorithm are able to converge to the global solutions with a relative higher speed and higher possibility. The simulations and experiments have demonstrated the effectiveness of our SVM/LM combined method.

## ACKNOWLEDGMENTS

This work was funded by the National Natural Science Foundation of China (51475191, 51405172, 51575214, and 51525502), the Natural Science Foundation of Hubei Province of China (2015CFB278 and 2015CFA005), the China Postdoctoral Science Foundation (2014M560607 and 2015T80791), and the Program for Changjiang Scholars and Innovative Research Team in University of China (Grant No. IRT13017).

**Tables**

Table 1. The extracted results of the 1D Si grating together with the time consumptions by the SVM/LM method, the LM algorithm, and SEM, respectively.

| Method | Initialization (nm) | TCD (nm) | Hgt (nm) | BCD (nm) | Time of SVM mapping (s) | Time of LM iterations (s) | Number of LM iterations |
|---|---|---|---|---|---|---|---|
| SVM/LM | (362.5, 487.5, 362.5) | 345.7 | 472.6 | 395.1 | 0.49 | 122.01 | 12 |
| LM | (400, 450, 400) | 345.7 | 472.6 | 395.1 | —— | 142.4 | 14 |
|  | (520, 430, 450) | 488.5 | 428.1 | 449.7 | —— | 1123.2 | 113 |
|  | (500, 525, 500) | 533.9 | 469.6 | 597.9 | —— | 513.8 | 52 |
|  | (475, 525, 500) | 407.7 | 512.4 | 665.8 | —— | 775.8 | 76 |
| SEM | —— | 350 | 472 | 383 | —— | —— | —— |

Table 2. The extracted results of the deep-etched multilayer grating together with the time consumptions by the SVM/LM method, the LM algorithm, and TEM, respectively.

| Method | Initialization (nm) | $D_1$ (nm) | $H_1$ (nm) | $D_2$ (nm) | $H_2$ (nm) | $D_3$ (nm) | $H_3$ (nm) | Time (s) |
|---|---|---|---|---|---|---|---|---|
| SVM/LM | (80, 130, 91, 16.75, 129, 129) | 75.86 | 128.85 | 90.37 | 16.49 | 125.20 | 131.76 | 354.78 |
| LM | (55, 115, 66, 1, 104, 114) | 91.95 | 77.21 | 105.50 | 51.95 | 120.98 | 104.83 | $2.33 \times 10^4$ |
|  | (95, 155, 106, 19, 144, 154) | 100.00 | 153.87 | 93.04 | 20.37 | 142.38 | 155.01 | $2.75 \times 10^4$ |
| TEM | —— | 75.01 | 135.60 | 86.90 | 9.92 | 124.13 | 134.29 | —— |

**Figures**

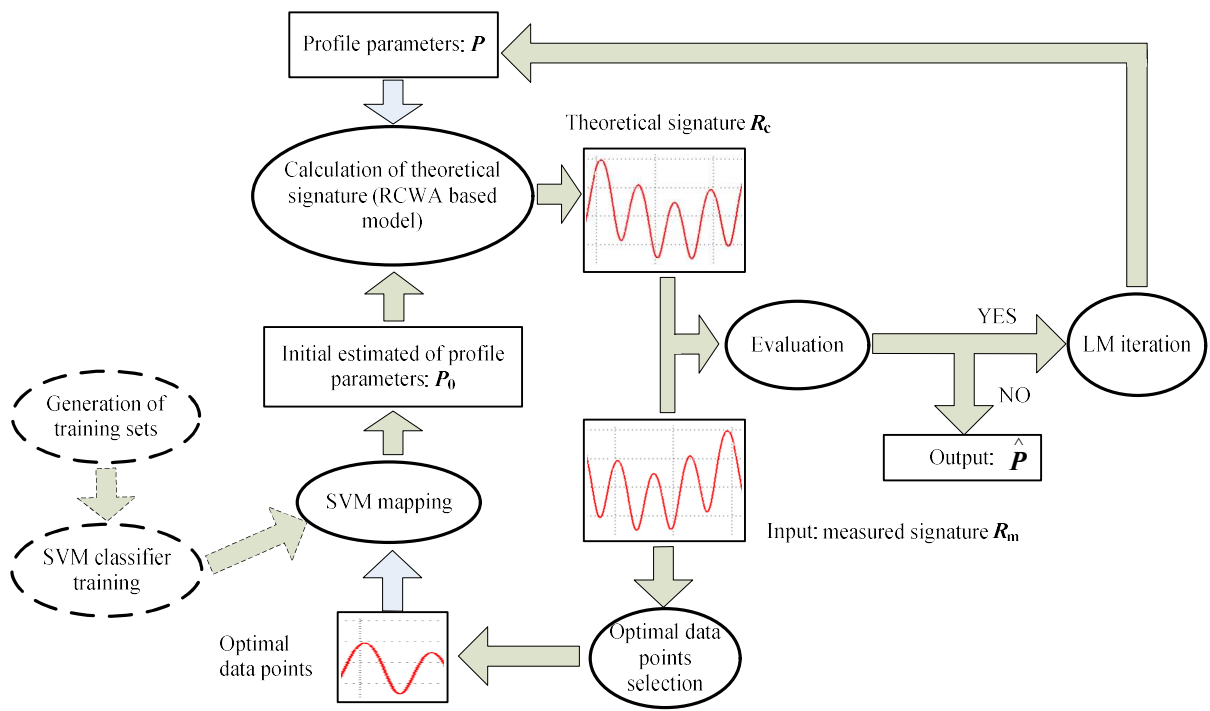

Fig. 1 (Color online) Flowchart of profile parameter extraction using the SVM/LM combined method.

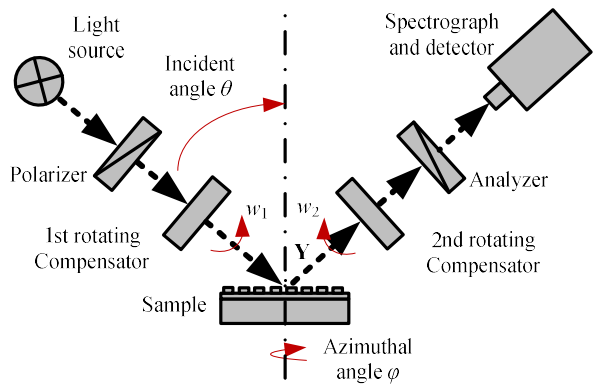

Fig. 2 (Color online) Principle of the dual rotating-compensator MME.

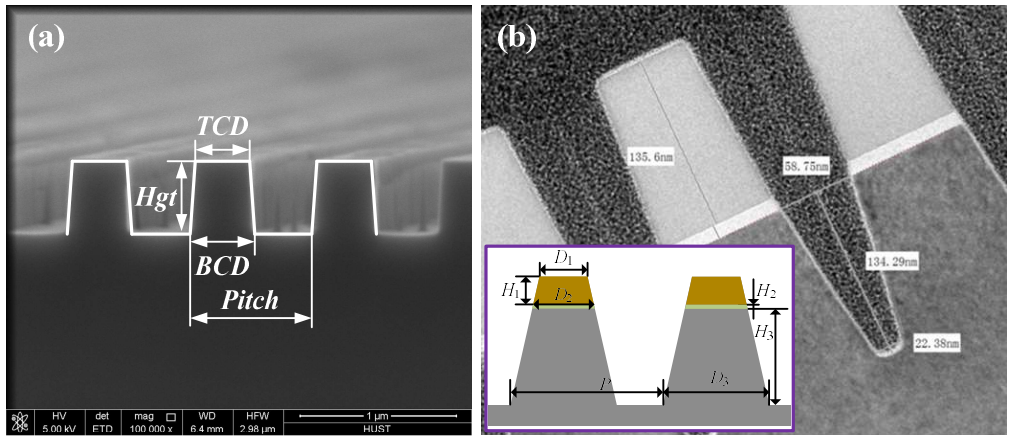

Fig. 3 Cross-section images of (a) the Si only grating and (b) the deep-etched multilayer grating consisting of Si, SiO$_2$, and Si$_3$N$_4$ trapezoidal gratings from bottom to top.

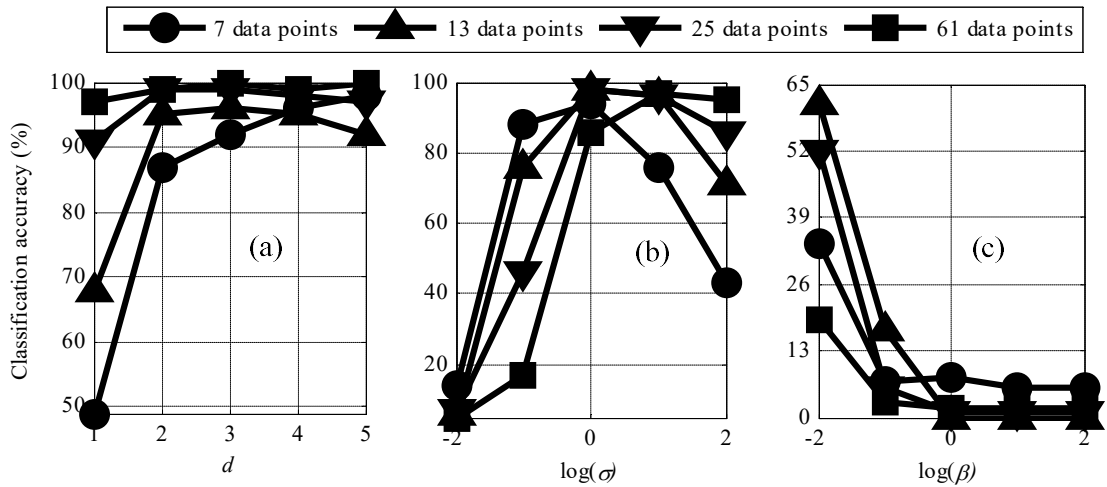

Fig. 4 Classification accuracy for the 1D Si grating as function of (a) $d$, (b) $\log(\sigma)$, and (c) $\log(\beta)$ in the polynomial, RBF, and Sigmoid kernels respectively. The solid lines with circles, upper triangles, lower triangles, and squares represent the 7, 13, 25, and 61 wavelength points in a training signature respectively. The number of training pairs correspond to each class is 3375.

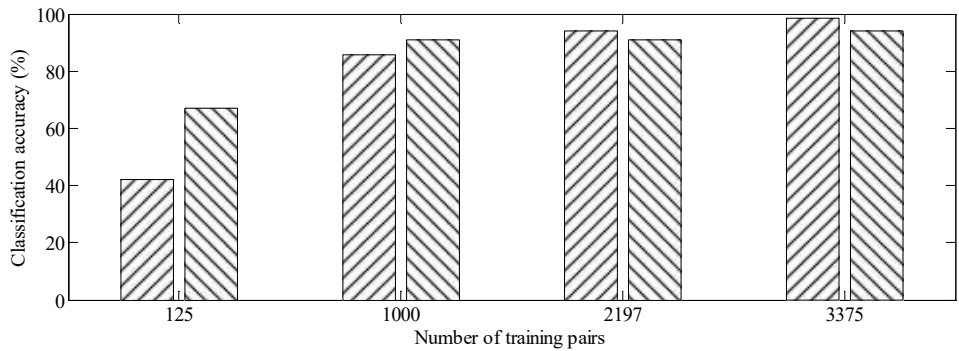

Fig. 5 (Color online) Classification accuracy for the 1D Si grating under different numbers of training pairs. The left and right slash marks represent the classification results obtained by the polynomial kernel with $d = 5$ and by the RBF kernel with $\sigma = 1$, respectively.

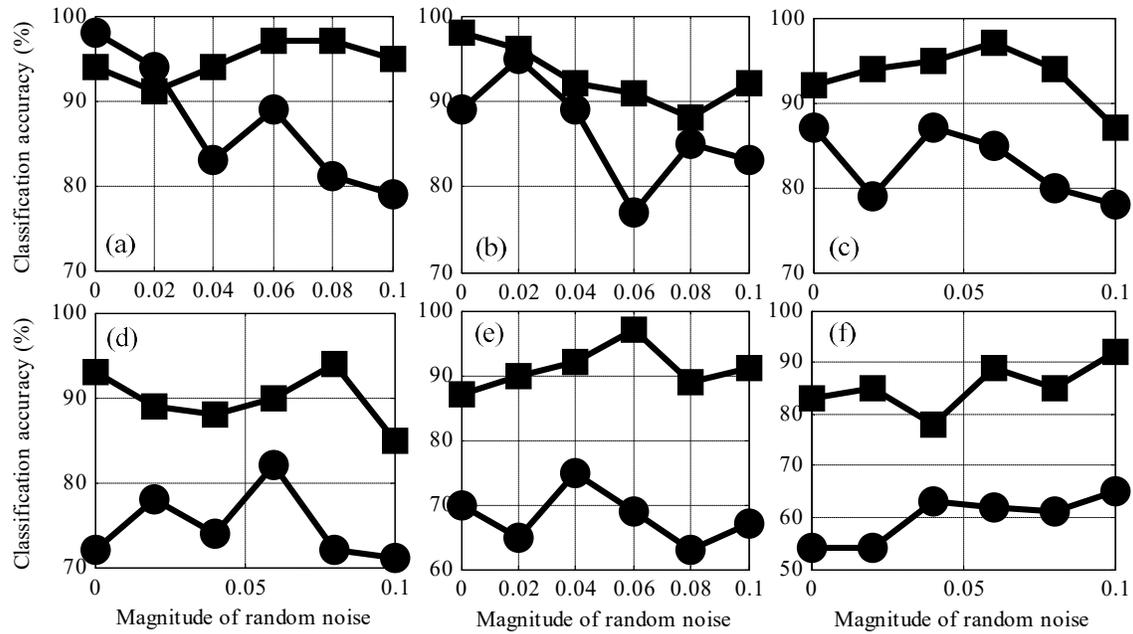

Fig. 6 Classification accuracy for the 1D Si grating as function of magnitudes of random errors under different offset errors, whose magnitudes are (a) 0, (b) 0.02, (c) 0.04, (d) 0.06, (e) 0.08, and (f) 0.10. The curves marked by circles and squares represent the results by the polynomial kernel with $d = 5$ and by the RBF kernel with $\sigma = 1$, respectively.

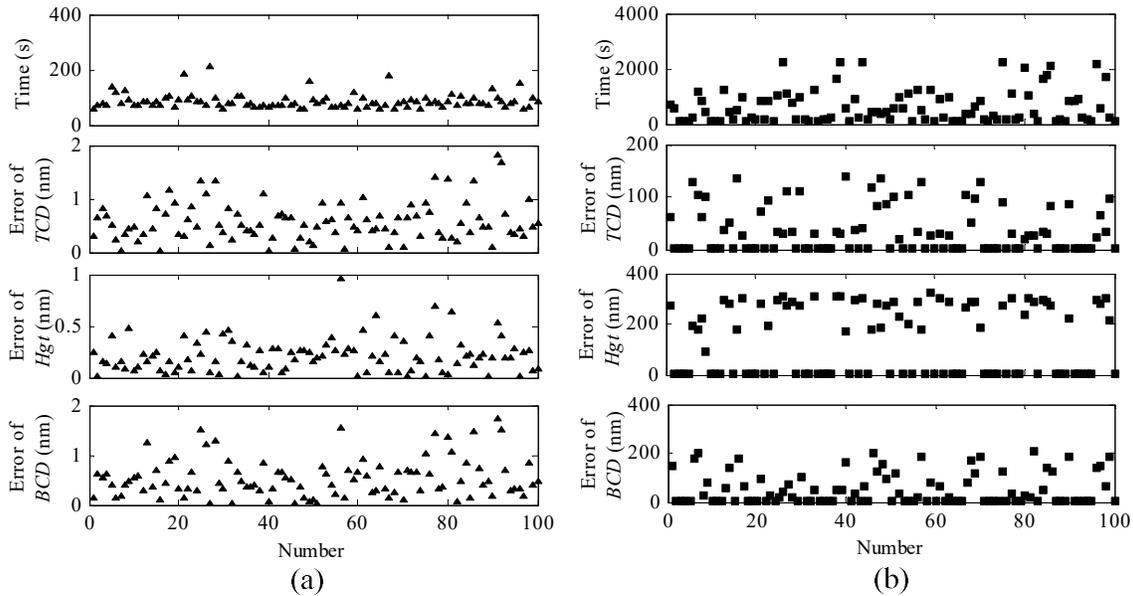

Fig. 7 Time consumptions and extracted errors of the 1D Si grating by (a) the proposed SVM/LM algorithm and (b) the LM algorithm respectively. The subplots from top to bottom correspond to the time consumption, extracted errors of $TCD$, $Hgt$ and $BCD$ respectively.

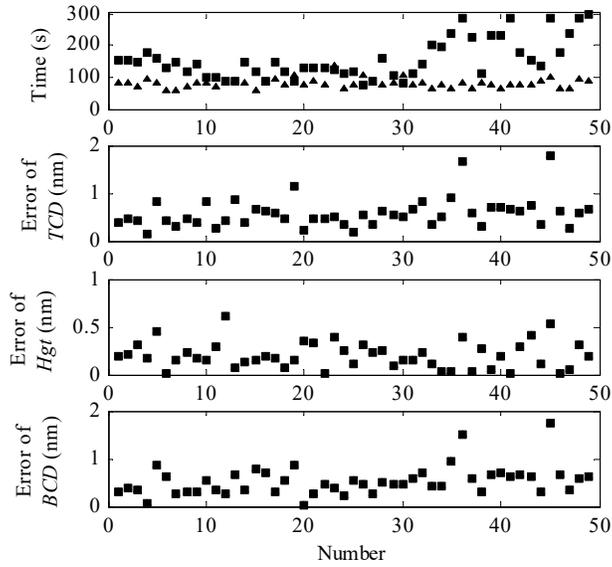

Fig. 8 Time consumptions and extracted errors of *TCD*, *Hgt* and *BCD* by SVM/LM and LM methods respectively.

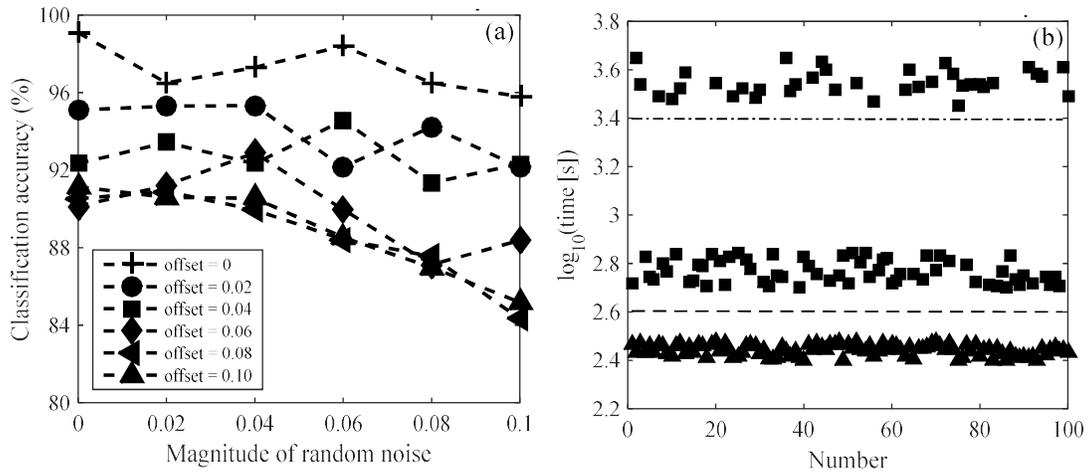

Fig. 9 (a) Classification accuracy for the deep-etched multilayer grating as function of different levels of random errors and offset errors; (b) Time consumptions of SVM/LM and LM iterations, respectively. The upwards triangular and square respectively represent the time consumptions of SVM/LM and LM methods.

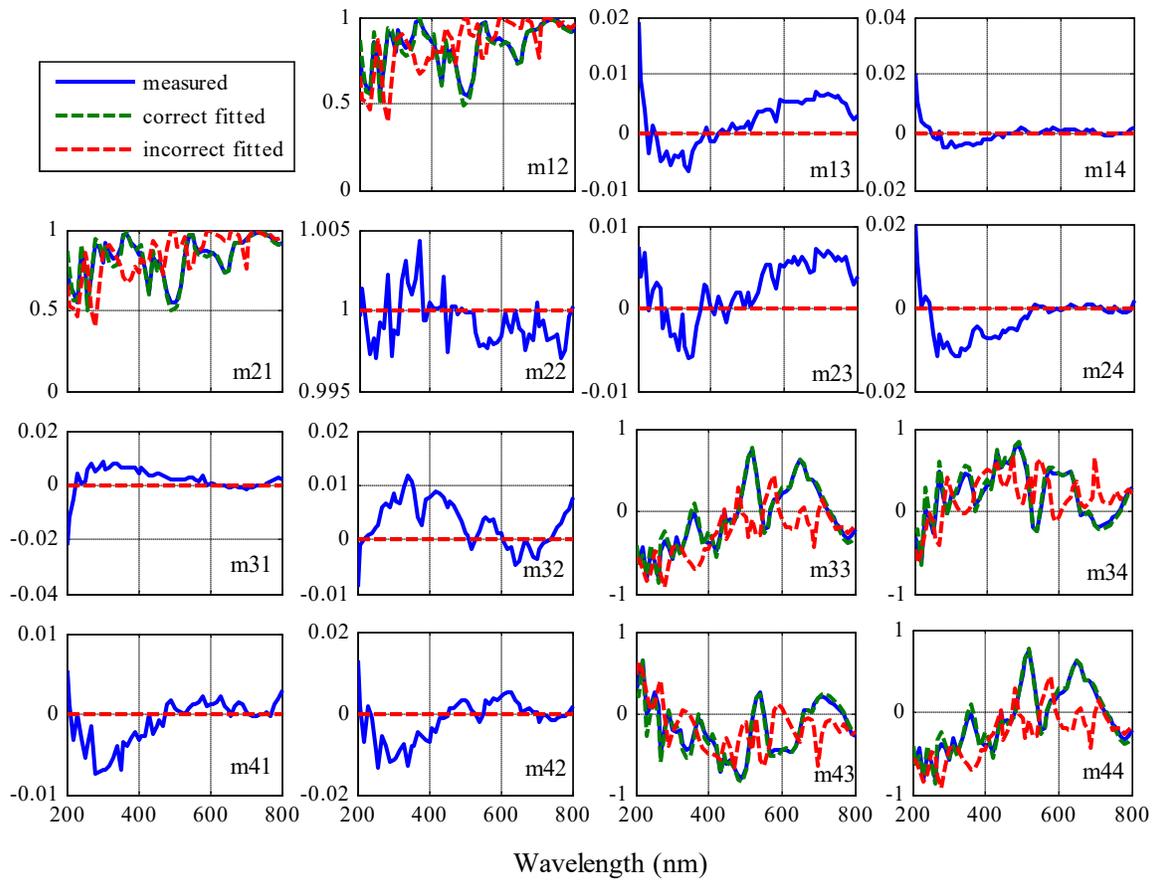

Fig. 10 (Color online) Comparison of the correctly and incorrectly fitted Mueller matrices of the measured Mueller matrix. The correctly fitted Mueller matrix represented by the green dotted line is obtained at an initialization after the SVM classifiers mapping, while the incorrectly fitted Mueller matrix represented by the red dash dotted line is obtained at the initialization (520, 430, 450) nm.